\documentclass{amsart}
\usepackage{amsfonts, amsmath, latexsym, epsfig}
\usepackage[norelsize]{algorithm2e}
\usepackage{amssymb}
\usepackage{color}
\usepackage{epsf}
\usepackage{url}

\newtheorem{proposition}{Proposition}
\newtheorem{theorem}{Theorem}

\newtheorem{conjecture}{Conjecture}

\def\QuotS#1#2{\leavevmode\kern-.0em\raise.2ex\hbox{$#1$}\kern-.1em/\kern-.1em\lower.25ex\hbox{$#2$}}

\begin{document}

\author{Adel Alahmadi}
\address{Math. Dept., King Abdulaziz University, Jeddah 21589, Saudi Arabia}
\email{adelnife2@yahoo.com}

\author{Husain Alhazmi}
\address{Math. Dept., King Abdulaziz University, Jeddah 21589, Saudi Arabia}
\email{alhazmih@yahoo.com}

\author{Shakir Ali}
\address{Math. Dept, Faculty of Science, Rabigh Campus, King Abdulaziz University, Rabigh, Saudi Arabia}
\email{shakir50@rediffmail.com}

\author{Michel Deza}
\address{Michel Deza, \'Ecole Normale Sup\'erieure, 75005 Paris, France}
\email{Michel.Deza@ens.fr}

\author{Mathieu Dutour Sikiri\'c}
\address{Mathieu Dutour Sikiri\'c, Rudjer Boskovi\'c Institute, Bijenicka 54, 10000 Zagreb, Croatia, Fax: +385-1-468-0245}
\email{mathieu.dutour@gmail.com}

\author{Patrick Sol\'{e}}
\address{Patrick Sol\'e, T\'el\'ecom Paris Tech, 46 Rue Barrault, 75013 Paris, France}
\email{patrick.sole@telecom-paristech.fr}

\title{
Hypercube emulation of  interconnection
networks topologies}
\date{}

\maketitle

\begin{abstract}
We address various topologies (de Bruijn, chordal ring, generalized Petersen,
meshes) in various ways ( isometric embedding, embedding up to scale, embedding up
to a distance) in a hypercube
or a half-hypercube. Example of obtained embeddings:  infinite series of hypercube  embeddable  Bubble Sort and Double Chordal Rings topologies, as well as of regular maps.
\end{abstract}

\section{Introduction}
The hypercube topology is a very popular topology for Parallel Processing computers
from the Connection Machine \cite{Kotsis} onward.
One way to emulate an alternative topology on such a computer is to address the
vertices of the guest topology by the vertex of the host hypercube or some subgraph
thereof.
This addressing can be used for routing purposes, for instance. 
Another important application is the addressing of knowledge databases \cite{SGL}.
This latter application is important for natural language processing.

In the present work we model graph theoretically the addressing process in various
ways from isometric embedding (the guest graph is a so-called {\em partial cube}) to embedding
up to scale
(geodetic distance on the host is a constant times that of the guest) or up to a given
distance (called henceforth
{\em truncated embedding}). This work is an application, a continuation  and a generalization of the book
\cite{DGS}, which considers only embeddings.
We shall consider many  popular topologies in turn and will question their
embeddability. The material is organised as follows.
To begin with, we consider  insertion/deletion-based distances in Section \ref{Indel} and other graphs defined on alphabets 
(Odd graph, Generalized Petersen, De Bruijn)
 in Section \ref{Alp}
and move on to cycle-based topologies in Section \ref{Cyc}.
Hypercube based topologies (Cube-connected Cycles, Butterfly graphs)  in Section \ref{Hyp}
and Cayley graphs on the group of permutations in Section \ref{Cay} are also
considered. 
We conclude in Section \ref{Map} by regular maps: 
skeletons of 
Klein graph, Dyck graph and so on.

\section{Preliminaries}

Denote by $H_m$ the skeleton of the $m$-dimensional cube. It is the graph on all binary sequences of length $m$
with two of them, say, $x=(x_1, \dots , x_n)$ and $y=(y_1, \dots , y_n)$, 
 being adjacent if their {\em Hamming distance} $$d_H(x,y)=\sum_{i=1}^m|x_i-y_i|$$ is $1$. Clearly, Hamming distance
 is an $l_1$-metric and the square of $l_2$-metric (Euclidean distance) on these sequences.

Denote by $\frac{1}{2}H_m$
the {\em $m$-half-cube graph}. It is defined on all  binary sequences of length $m$ having even number of ones,
with two of them, say, $x$ and $y$, 
 being adjacent if 
 $d_H(x,y)=2$. Clearly, $H_m$ is an isometric subgraph of 
$\frac{1}{2}H_{2m}$. 

Given a finite connected graph $G=(V,E)$ of diameter $d$, let $D(G)$ denote the (shortest path) distance matrix
$((d_{ij}))$ of its vertices. Call $G$ {\em embeddable} (or, as it is done in \cite{BCN},
 {\em code graph})  and denote it by $G\to \frac{1}{2}H_m$,
if $G$ is an isometric subgraph of a some $m$-half-cube, i.e., $((d_{ij}))$ embeds {\em scale-$2$-isometrically} into the distance matrix of  $m$-cube. 
If, moreover, $G$ is an isometric subgraph of a  $m$-cube,
 denote it by $G\to H_m$ and call $G$ a {\em partial cube}.
Clearly, $$\mbox{~if~}\,\, G\to \frac{1}{2}H_m\,\, \mbox{~and~}\,\,\, G'\to \frac{1}{2}H_{m'},\,\,\,\mbox{~then~}\,\,\, G\times G' \to \frac{1}{2}H_{m+m'}.$$
Another isometric subgraph of $\frac{1}{2}H_{m}$ is the {\em Johnson graph} $J(m,k)$; its vertices  are the $k$-element subsets of an $m$-element set, and two vertices are adjacent when they meet in a $(k-1)$-element set.
Let us denote by $G\to J(m,k)$ such eventual special case of embedding into $\frac{1}{2}H_{m}$.

\begin{theorem}(\cite{tylkin})  

For  a connected graph $G$, it holds:

(i) $G$ is {\em $l_1$-embeddable} (i.e., it embeds isometrically into some $l_1$-space) 
 if and only if $D(G)$,
 for some integers $m, \lambda \ge 1$,  embeds
  {\em scale-$\lambda$-isometrically}  into the distance matrix of  $m$-cube;

(ii) if $G$ is  $l_1$-embeddable, then it is  an {\em hypermetric graph}, i.e., 
its  $d_G$ satisfies all {\em hypermetric inequalities}
$$\sum_{1\leq i<j\leq n}b_i b_j \: d_G(v_i,v_j)\leq 0,$$
where $b=\{b_1,b_2,\dots,b_n\}\in Z^n$, $\sum_{i=1}^n b_i=1$ and $v_i$,
$1\le i\le n$, are vertices of $G$.
\end{theorem}

\begin{figure}
\begin{center}
\setlength{\unitlength}{1cm}
\begin{minipage}[t]{5.5cm}
\hfil\begin{picture}(5,5)
\leavevmode
\epsfxsize=5.5cm
\epsffile{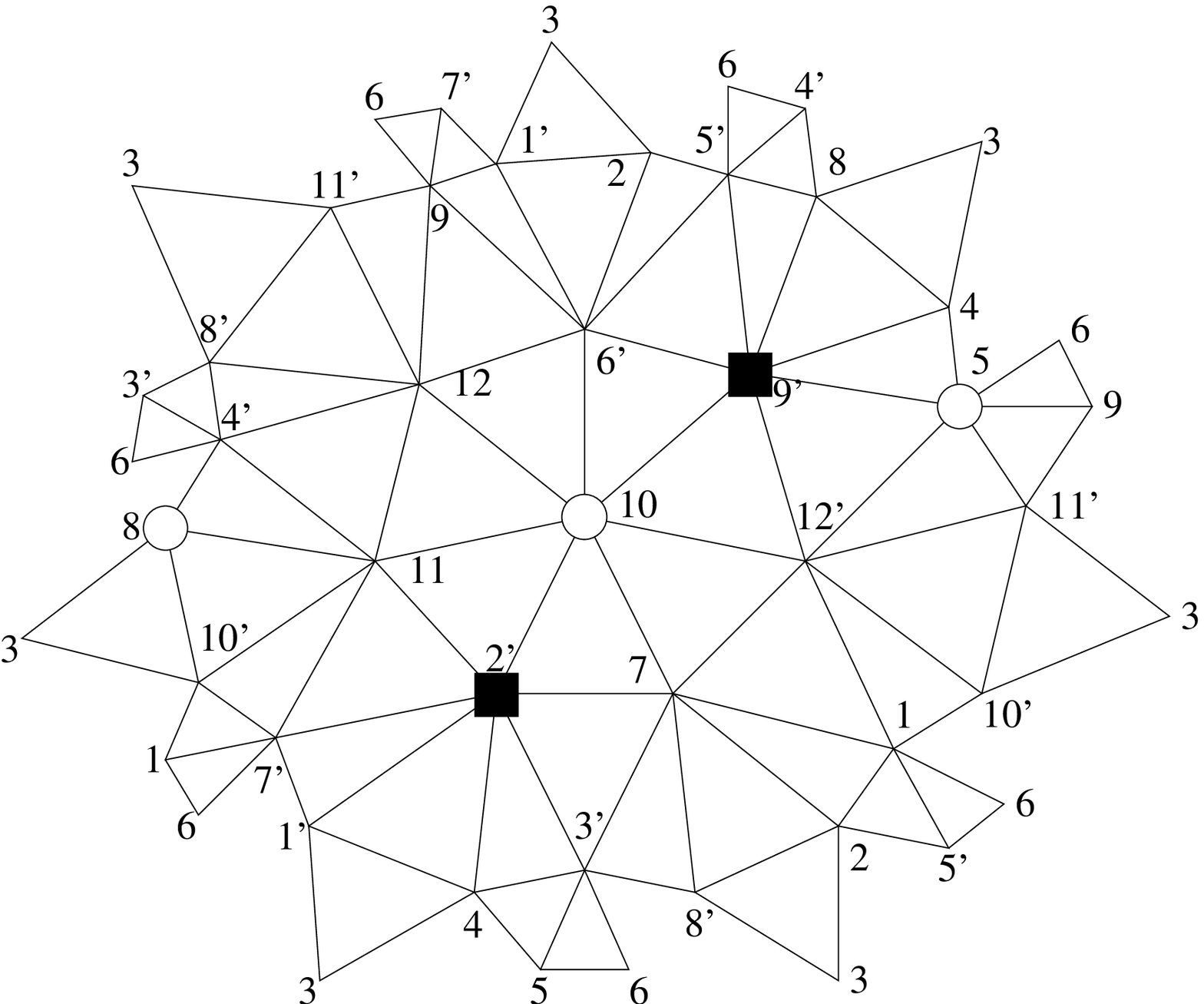}
\end{picture}
\hfil\par
\end{minipage}
\hspace{1mm}
\begin{minipage}[t]{4.3cm}
\hfil\begin{picture}(3,3)
\leavevmode
\epsfxsize=43mm
\epsffile{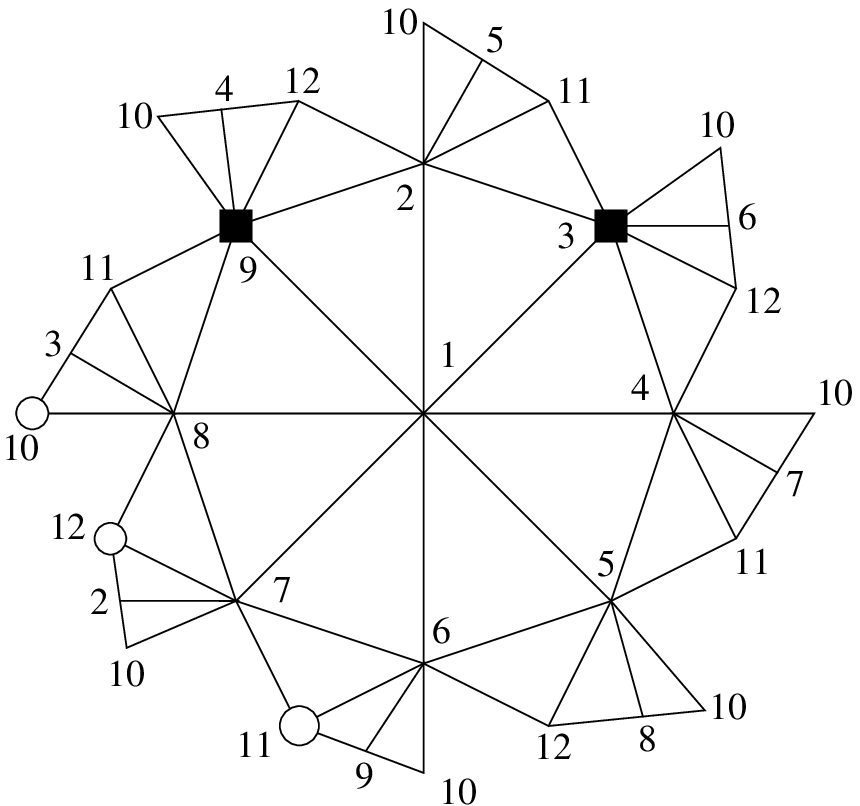}
\end{picture}\hfil\par
\end{minipage}
\end{center}

\caption{Examples of not $5$-gonal regular maps: dual Klein map $\{7,3\}$ and dual Dyck map $\{8,3\}\simeq K_{4,4,4}$
 on genus $3$ surface}
\label{Dyck}
\end{figure}

The  inequality with $\sum_{i=1}^n|b_i|=2k+1$
is called a {\em $(2k+1)$-gonal
inequality}. 
Clearly, the case $k=1$ corresponds to the usual triangle inequality. 
The  {\em 5-gonal inequality}
correspond to  $b_a=b_b=b_c=1$, $b_x=b_y=-1$, i.e., it is 
$$d(x,y) + (d(a,b)+d(a,c)+d(b,c)) \le \sum_{i=a,b,c}(d(x,i)+d(y,i)$$
 for 
 any vertices $a,b,c,x,y$.
 \cite {Avis} showed that a connected graph  is a partial cube if and only if it is  bipartite and its path-metric satisfy all $5$-gonal inequalities.
 See examples of not $5$-gonal graphs on Fig. \ref{Dyck}.

\begin{figure}\begin{center}
\begin{minipage}[t]{4cm}
\centering
\epsfig{file=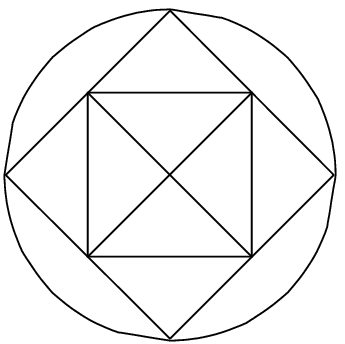}\par
\end{minipage}
\begin{minipage}[t]{4cm}
\centering
\epsfig{file=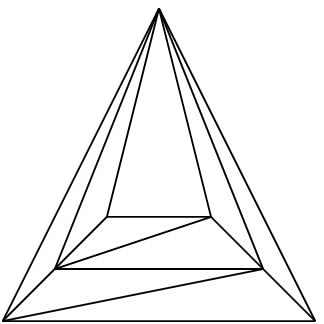}\par
\end{minipage}
\end{center}
\caption{Examples of hypermetric, but not embeddable graphs}
\label{fig:EmbeddingHyp}
\end{figure}

 The hypermetricity is not sufficient, if the number of vertices is greater than $6$, for embeddability and, larger, for $l_1$-embeddability; see Fig.
 \ref{fig:EmbeddingHyp}. 

\begin{theorem}(Theorem 17.1.1 in \cite{DL})
For  a connected graph $G$, it holds:

(i) $G$ is hypermetric if and only if it is an isometric subgraph of a Cartesian product of half-cube graphs $\frac{1}{2}H_m$, {\em cocktail-party graphs} $K_{2, \dots ,2}$ and copies of the {\em Gosset graph} $G_{56}$;

(ii) $G$ is an $l_I$-graph if and only if it is an isometric subgraph of  a Cartesian product of  half-cube graphs $\frac{1}{2}H_m$
 and cocktail-party graphs $K_{2, \dots ,2}$.
\end{theorem}

Given an integer $2\le s\le d$, call $G$ {\em $s$-tr.embeddable} (short for {\em up to $s$ truncated-embeddable}) if there exists a  distance matrix  $D'=((d'_{ij}))$ of order $|V|$ with $d'_{ij}=d_{ij}$, whenever $d_{ij}\le s$, which
 is isometrically embeddable in the distance matrix of some some $m$-half-cube. So, $D'$ is a graphic distance matrix
 only if $D'=D$.

 Clearly, 
 $0\le d'_{ij}\le d_{ij}$, whenever $d_{ij}>s$, and   $s$-tr.embeddability implies  $(s-1)$-tr.embeddability. 
Such embedding produces an addressing of vertices by binary sequences, preserving, up to scale $2$  and up to value $s$
 the graph's distances. This addressing can be used for routing purposes and
 for the emulation of an architecture of topology $G$ on a host machine
that has a  hypercube topology, as, say, 
 the Connection Machines from CM-1 to CM-4.

The classical structures, used as topologies for 
interconnection networks, are, say,  
trees, hypercubes and rings. First two, as well as even rings, are, trivially, partial cubes.
An odd ring $C_n$ embeds into $\frac{1}{2}H_n$. Here we consider eventual embedding for other topologies, especially, for those mentioned in 
a good surveys \cite{Kotsis}, \cite{He}. This work, while being a follow-up of the book \cite{DGS}, since we again investigate embeddability, is focused now on network applications.
Also, in the case of absence of embedding, we look now for "embedding" in weaker sense, i.e., maximal $s < d$, for which $s$-tr.embedding eventually exist.
 We not consider the setting of $l_1$-embedding, because it looks as not suitable for applications.

\section{Algorithm}
Our work is based on heavy computations, using programs based on algorithm in \cite{DS}.
This algorithm (with time complexity $O(n^2 + nm)$ and space complexity $O(n^2)$) constructs
an embedding into $H_m$ up to scale $2$, if one exists.
The method has been extended to scale $2$ embeddings up to a given distance $s$.

For any edge $e=\{v, v'\}$ of $G$ and any such embedding $\phi$, the difference $\phi(v) - \phi(v')$
corresponds to a set $S_e$ of length $2$.
Lemma 4.1 in \cite{DS} allows in some cases to compute the size $i(e,e')$ of the intersection $S_e \cap S_{e'}$ for
two edges $e$ and $e'$. One of the conditions of applicability of the lemma is that vertices in $e$ and $e'$
are at distance at most $s$.

In the case of $s=diam(G)$, we can take a spanning tree $T$ of $G$ and compute $i(e,e')$
for all pairs of edges $e, e'$. If the function $i$ is negative, then
the graph is not embeddable. Otherwise, we can identify edges, such that $i(e,e')=2$,
and check that this defines an equivalence relation. Afterwards we define a graph $H$ on classes $\overline{e}$
and $\overline{e'}$ with two classes adjacent if $i(e,e')=1$.
We then check if the graph admits an inverse line graph by implementing the algorithm of \cite{Roussopoulos}.

In some cases the embedding is not unique; see, for example,  Tetrahedron in Figure \ref{fig:Embeddings}. All such cases of non-unique reversed line graph are classified in \cite{Roussopoulos}.

In the general case of $s$-tr.embedding, we may not have computed all the values 
$i(e,e')$. However, in the case of $s$-tr.embedding, one has $i(e,e') \in \{0,1,2\}$ and also  following consistency
relations:
\begin{enumerate}
\item If $i(e_1, e_2)=2$, then for any other edge $e'$,  it holds $i(e_1, e') = i(e_2, e')$.
\item If $e_1$, $e_2$, $e_3$ are three edges with $i(e_i, e_j)=1$ for $i\not=j$, then the edges
$e_1$, $e_2$ and $e_3$ can be of the form:
\begin{equation*}
e_1=AB, e_2=AC, e_3=BC.
\end{equation*}
In that case, for any other edge $e$, we will have, up to permutation, following patterns of intersections: 
$\{i(e_1, e), i(e_2, e), i(e_3, e)\}$ : $\{1,1,0\}$, $\{1,1,2\}$ and $\{0,0,0\}$.
Alternatively, the edges $e_1$, $e_2$ and $e_3$ can be of the form:
\begin{equation*}
e_1=AB, e_2=AC, e_3=AD.
\end{equation*}
In that case, for any other edge $e$, there is, up to permutation, following patterns of intersections: 
$\{i(e_1, e), i(e_2, e), i(e_3, e)\}$ : $\{1,0,0\}$, $\{1,1,0\}$, $\{1,1,1\}$, $\{1,1,2\}$
and $\{0,0,0\}$.
Also, the number of patterns $\{1,1,0\}$ is at most $3$.
\end{enumerate}

At start, when Lemma 4.1 of \cite{DS} can be applied, we set up $i(e,e')$. If it cannot be applied, we only
know that $i(e,e') \in \{0,1,2\}$. With the above relations, one can sometimes deduce $i(e,e')$ from 
what is known and this can be iterated.
Therefore, in some cases the values of $i(e,e')$ is completely determined from the values obtained
from Lemma 4.1 in \cite{DS}.
But in other cases, the above logical relations are not sufficient to deduce all possible values
$i(e,e')$. Thus, we apply a classical backtracking strategy of choosing the value $i(e,e')$, applying
above deduction rules and iterating until we find all possible embeddings.

\section{Indel-based graphs}\label{Indel}

 Denote by $D_{n}$ the set of binary sequences of length $n$ and, for any $0\le i\le n$, denote by
$D_{i,\dots ,n}$  the set $\cup _{j=i}^{n}D(j)$.

The {\em indel  graph} $Ind_{i,\dots ,n}$ is defined on 
on $D_{i,\dots ,n}$ by considering two sequences adjacent if one can be obtained from the other by {\em indels}, i.e., insertions or deletions of characters only. This graph  is bipartite and has diameter $2n$.
\begin{proposition}
With exception of $Ind_{0,1}=P_2\to H_{2}$ and
$Ind_{n-1,n}=H_{2n}$ for $n=2,3$,
any $Ind_{i,\dots ,n}$ 
 is not $2$-tr.embeddable.
\end{proposition}
\proof In fact, see three embeddings  on Table \ref{tab1}. Clearly,  $Ind_{n-1,n}$
 is an isometric subgraph of     $Ind_{n,n+1}$ and, for $0\le i\le n-2$, of $Ind_{i,\dots ,n}$. 
By computation, $Ind_{0,1,2}$ and $Ind_{3,4}$  are not $2$-tr.embeddable, proving the Proposition. \qed

\begin{table}
\begin{center}
\begin{minipage}{6cm}
\begin{center}
{\small
\begin{tabular}{||c|c||}
\hline
\hline
vertex& vertex address\\
\hline
\hline
{\bf 0,0}      & ({\bf 0,0,0,0},0,0)\\
{\bf 1,0}      & ({\bf 0,1,1,0},0,0)\\
1,1            & (0,1,1,1,1,0)\\
0,1            & (0,0,1,0,1,0)\\
{\bf 0,0,0}    & ({\bf 1,0,0,0},0,0)\\
{\bf 1,0,0}    & ({\bf 0,1,0,0},0,0)\\
{\bf 0,1,0}    & ({\bf 0,0,1,0},0,0)\\
{\bf 1,1,0}    & ({\bf 0,1,1,1},0,0)\\
1,1,1          & (0,0,1,1,1,1)\\
1,0,1          & (0,1,1,0,1,0)\\
0,1,1          & (0,0,1,1,1,0)\\
0,0,1          & (0,0,0,0,1,0)\\
\hline
\hline
\end{tabular}
}
\end{center}
\end{minipage}
\begin{minipage}{6cm}
\begin{center}
\centering
\resizebox{35mm}{!}{\includegraphics{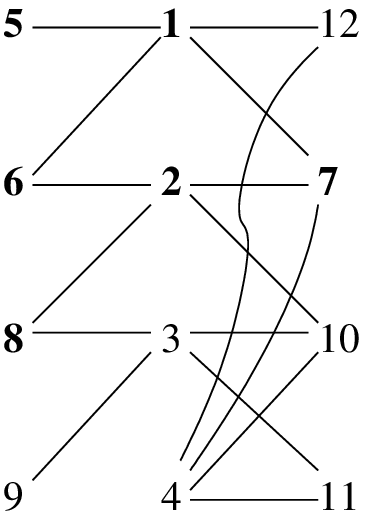}}\par
\end{center}
\end{minipage}

\end{center}
\caption{Embedding  $Ind_{2,3}\to H_6$. The  boldfaced minor of
 rows $1,2,5-8$ and columns $1-4$ gives
embedding  $Ind_{1,2}\to H_4$. The minor of
 rows $1,2,5$  and columns $1,2$ gives
embedding $Ind_{0,1}\to H_2$}

\label{tab1}
\end{table}

The {\em Levenshtein  graph} $Lev_{i,\dots ,n}$ is defined on 
on $D_{i,\dots ,n}$ by considering two sequences adjacent if one can be obtained from the other by  indels and changes of characters as  say, $x$ on $1-x$, only. This graph  has diameter $n$.
\begin{proposition}
With exception of $Lev_{0,1}=K_3\to \frac{1}{2}H_{3}$,
any $Lev_{i,\dots ,n}$ 
 is not $2$-tr.embeddable.
\end{proposition}
\proof In fact,   $Lev_{n-1,n}$
 is an isometric subgraph of     $Lev_{n,n+1}$ and, for $0\le i\le n-2$, of $Lev_{i,\dots ,n}$. 
By computation, $Lev_{1,2}$ is not $2$-tr.embeddable, proving the Proposition. \qed

Note, that \cite{Andoni} gave lower bound $\frac{3}{2}$ for distortion of  a 
{\em Lipchitz $l_1$-embedding} of Levenshtein metric on sequences; see also Problem 2. 15 in \cite{Mat}.
Given metric spaces $(X,d_X)$ and $(Y,d_Y)$, the {\em distortion} of  a mapping $f:X \rightarrow Y$  is $||f||_{Lip}||f^{-1}||_{Lip}$,
where the {\em Lipschitz norm} is defined by
$$f_{Lip}=\sup_{x,y\in X, x\neq y} \frac{d_Y(f(x),f(y))}{d_X(x,y)}.$$

The {\em Ulam metric} (or {\em permutation editing metric}) $U$
 is an editing metric on $Sym_n$, obtained by 
{\em character moves}, i.e.,  transpositions of characters.
 It is the half of the indel metric on $Sym_n$ and It is right-invariant.
  Also, $n-U(x, y) = LCS(x, y)$, where $LCS(x, y)$ is the length of the longest common subsequence (not necessarily a substring) of $x$ and $y$.
The {\em Ulam graph} $Ul(n)$ has diameter $n-1$.

\begin{proposition}
With exception of  $Ul(2)=K_2=H_1$ and $Ul(3)=K_{2,2,2}\to \frac{1}{2}H_4$,
any $Ul(n)$ is not $2$-tr.embeddable.
\end{proposition}
\proof In fact, see the embedding of $Ul(3)$, i.e., Octahedron,  in Fig. \ref{fig:Embeddings}
Clearly, $Ul(n)$ is an isometric subgraph of $Ul(n+1)$. 
By computation $Ul(4)$  is not $2$-tr.embeddable, proving the Proposition. \qed

\section{Network graphs on alphabets}
\label{Alp}

Here we consider some graphs, where the vertices are labeled by words of length $n$ over an alphabet and their relatives.

The {\em Odd graph} $O_n$ has one vertex for each of the $(n-1)$-element subsets of a $(2n-1)$-element set; two vertices are connected by an edge if and only if the corresponding subsets are disjoint. Any $O_n$
is a distance-transitive graph of diameter $n-1$. 

The Petersen graph is $O_3$. Any $O_n$ with $n\ge 4$ is even not $3$-tr.embeddable, since $O_4$ is not embeddable
and any $O_{n-1}$ is an isometric subgraph of $O_n$.
But for the bipartite double of $O_n$, called {\em Double Odd graph} (or {\em revolving doors}) $DO_{2n-1}$,
it holds
$DO_{2n-1}\to H_{2n-1}$;
together with hypercubes and even cycles, those graphs are only distance-regular ones (\cite{Koo}), which are partial cubes.
Note that $DO_5$ is the Desargues graph $GP(10,3)$.

The {\em Generalized Petersen graph} $GP(n,k)$
 is (Coxeter, 1950) 
 a graph consisting of an inner star polygon $\{n,k\}$and an outer regular polygon $\{n\}$ with corresponding vertices in the inner and outer polygons connected with edges. 
For example, $GP(5,2)$,
$GP(8,3)$ and $GP(12,5)$ are well-known  {\em Petersen graph}, {\em M\"obius--Kantor graph} and {\em Nauru graph}, respectively.

All case of embeddable $GP(n,k)$ are given in Table \ref{tab2}.
M\"obius--Kantor graph and Nauru graph are not embeddable and, moreover, not  $3$-tr.embeddable.

\begin{table}
\begin{center}
\begin{tabular}{||c|c|c||c||}
\hline
\hline
(n,k)&graph's name if any&diameter&embedding  into\\
\hline
(n=2m,1)&$Prism_{2m}$&m+1&$H_{m+1}$\\
(n=2m+1,1)&$Prism_{2m+1}$&m+1&$\frac{1}{2}H_{2m+3}$\\
(10,3)&Desargues graph&5&$H_5$\\
(10,2)&Dodecahedron&5&$\frac{1}{2}H_{10}$\\
(9,2)&&4&$\frac{1}{2}H_9$\\
(6,2)&D\"urer octahedron&4&$\frac{1}{2}H_8$\\
(5,2)&Petersen graph&2&$\frac{1}{2}H_6$\\
\hline
\hline
\end{tabular}
\end{center}
\caption{The cases of embedding for Generalised Petersen graph $GP$(n,k)}
\label{tab2}
\end{table}

\begin{figure}\label{DurerOctahedronElogatedDodecahedron}
\begin{center}
\begin{minipage}{4.1cm}
\centering
\epsfig{file=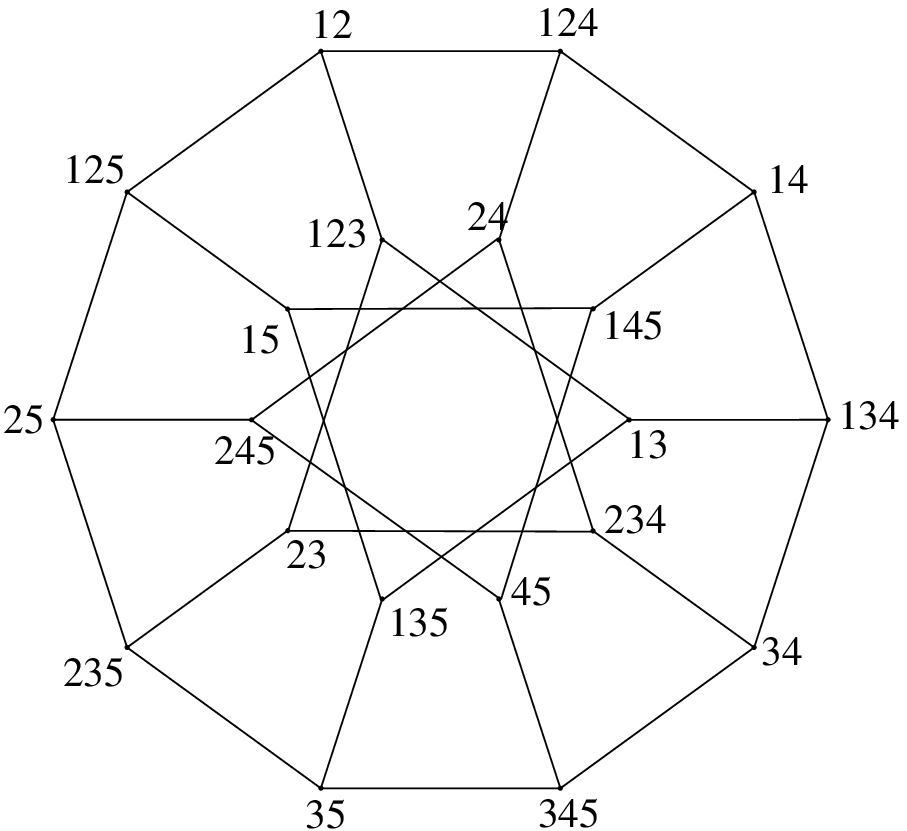,height=4.2cm}\par
\end{minipage}
\begin{minipage}{4.1cm}
\centering
\epsfig{file=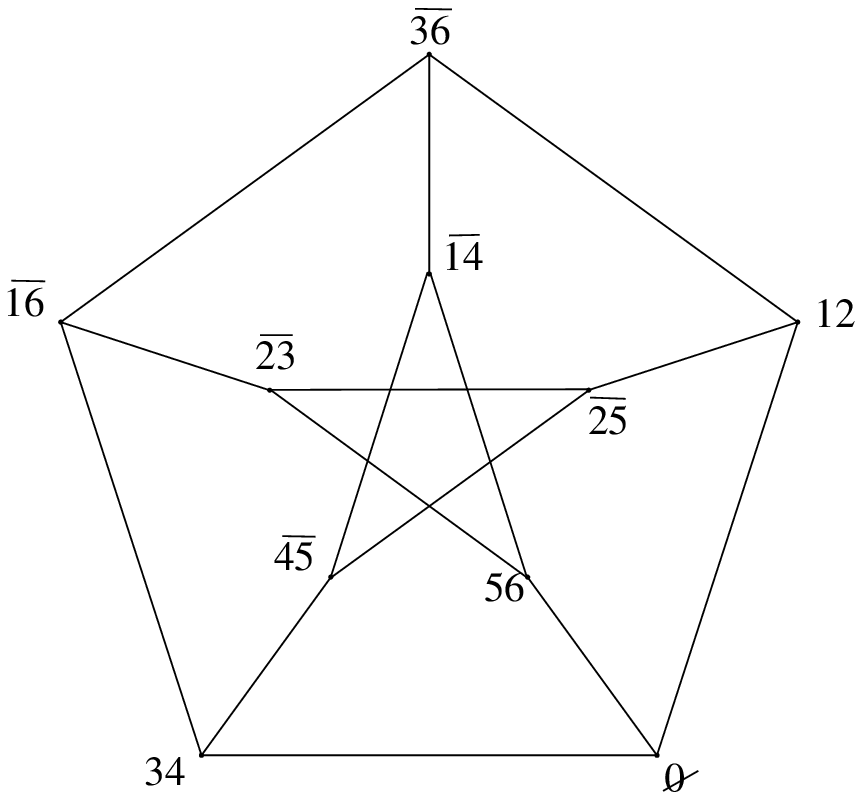,height=4.2cm}\par
\end{minipage}

\begin{minipage}{4.1cm}
\centering
\epsfig{file=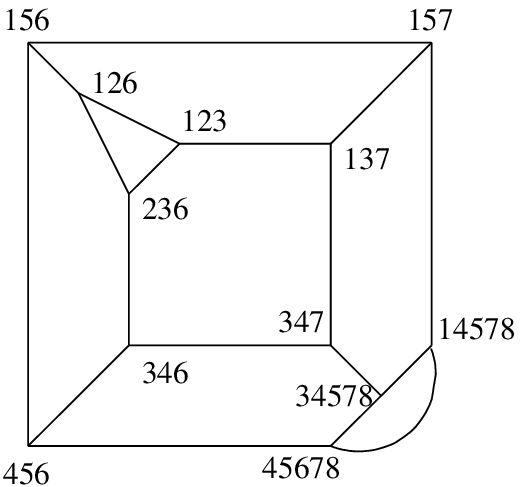,height=4.2cm}\par
\end{minipage}
\begin{minipage}{4.1cm}
\centering
\epsfig{file=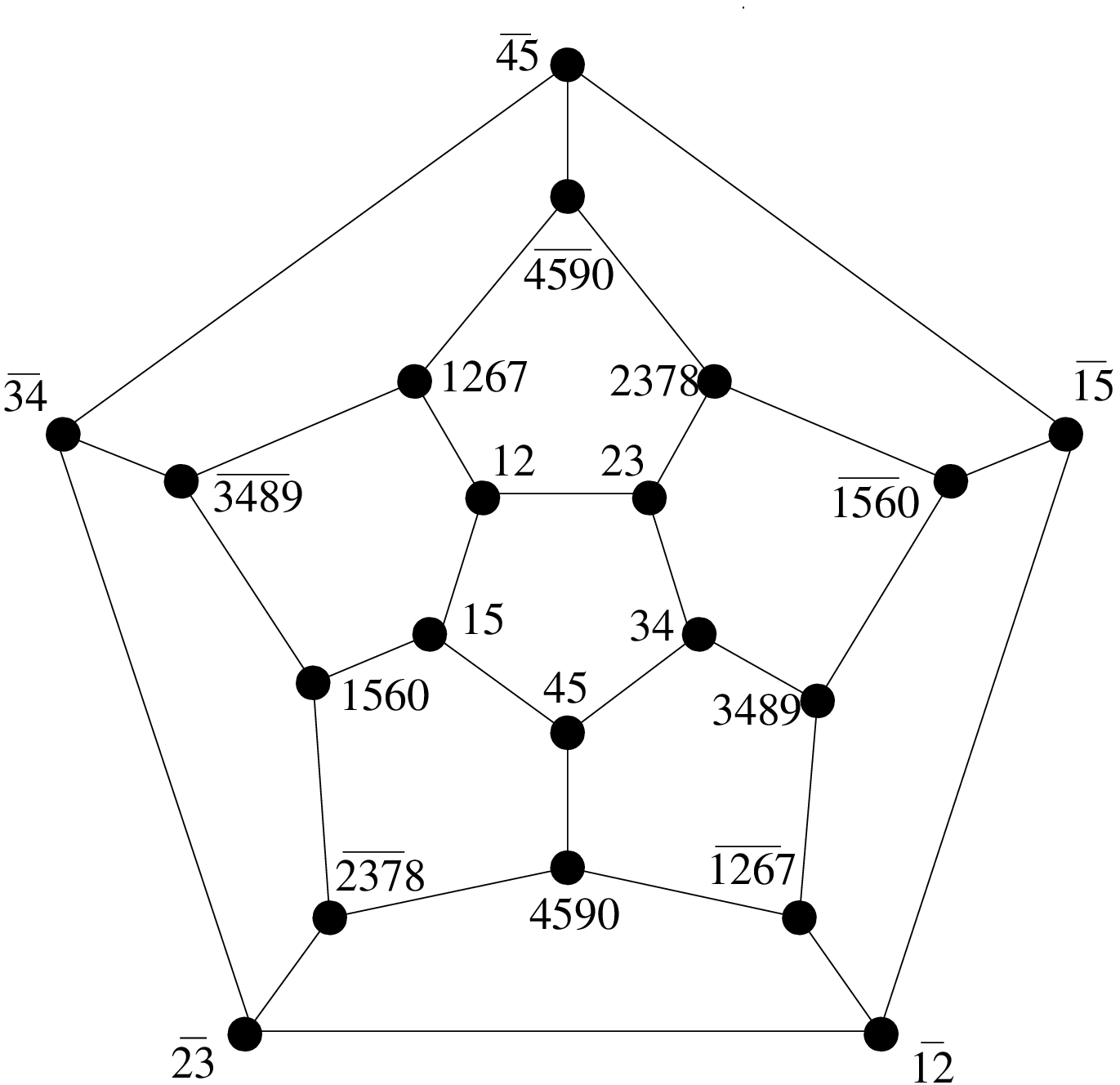,height=4.2cm}\par
\end{minipage}
\caption{Embedding of the Desargues Graph $GP(10,3)$ into $H_5$ and of the 
Petersen graph $GP(5,2)$, D\"urer's octahedron $GP(6,2)$
and Dodecahedron $GP(10,2)$
into $\frac{1}{2}H_6$, $\frac{1}{2}H_8$, $\frac{1}{2}H_{10}$,
 respectively}
\end{center}
\end{figure}

A {\em Moore graph} is a graph of diameter $d$ with girth $2d + 1$. The Moore graphs are:
 $K_n (n>2)$, $C_{2n+1}$,
the Petersen graph,
 the {\em Hoffman--Singleton graph} $HS$ (diameter $2$, girth $5$, degree $7$, order $50$) and 
a hypothetical graph of diameter $2$, girth $5$, degree $57$ and order $3,250$.
We found that $HS$ is not embeddable.

The undirected  {\em De Bruijn graph} $Br(m,n)$ is a 
graph on  $m^n$  $n$-tuples 
$(a_1 \dots a_n)$
over a $m$-character alphabet (denoted by juxtaposition). The edges are defined to be 
pairs of the form 
$((a_1 \dots a_n), (a_2 
 \dots a_na_{n+1}))$, where $a_{n+1}$ is any character in the alphabet. 
The undirected  {\em Kautz graph} $Ka(m,n)$ 
is defined similarly, but only tuples $(a_1 \dots a_n)$
with $a_i\neq a_{i+1}$ for each $i$ are taken. 
 The diameters of 
$Br(m,n)$ and, for $m\ge 3$, $Ka(m,n)$  are $n$. Clearly, $$Br(2,2)=K_4-e\to \frac{1}{2}H_4\,\,\,\,\,\mbox{~and~}\,\,\,\,\, Ka(2,n)=K_2=H_1.$$

\begin{conjecture}
\vspace{2mm}

(i) all  $Br(m,n)$ with $(m,n)\neq(2,2)$ are not $2$-tr.embeddable;  we checked it for
$(m,n)=(3,2)$, $(4,2)$, $(5,2)$, $(2,3)$, $(3,3)$, $(4,3)$, $(5,3)$,
$(2,4)$, $(3,4)$, $(4,4)$,$(2,5)$, $(3,5)$, $(3,6)$;

(ii) all $Ka(m,n)$ with $m\ge 3$ are not $2$-tr.embeddable; we checked it for
$(m,n)=(3,2)$, $(4,2)$, $(5,2)$, $(6,2)$, $(3,3)$, $(4,3)$, $(5,3)$,
$(3,4)$, $(4,4)$, $(5,4)$.
\end{conjecture}

\section{Hypercube Structures}\label{Hyp}

The {\em Cube-connected Cycles} $CCC_n$ is \cite{PrepVuil} cubic graph, formed by replacing
each vertex of an $n$-cube graph by a $n$-cycle. So, for example, $CCC_3$ is Truncated Cube.
The diameter of $CCC_n$ is $6$ for $n=3$ and $\left\lfloor \frac{5n-4}{2} \right\rfloor$ for $n\ge 4$;
note that value $\left\lfloor \frac{5n-2}{2} \right\rfloor$, given in Table 3.4
of \cite{Kotsis}, is correct only for $n=3$.  
\begin{conjecture} (checked $n=3,4,5$)

$CCC_n$ is not $(n+1)$-tr.embeddable.
\end{conjecture}
$CCC_3$ is the only planar $CCC_n$; its dual embeds into $\frac{1}{2}H_{12}$.
In fact (cf. \cite{DGS}), for any semiregular polyhedron $P$ (i.e.,  one of
$13$ Archimedean polyhedra, prisms and antiprisms),
exactly one of skeletons of $P$ and its spherical dual $P^*$ is embeddable.

A {\em Generalized Boolean $n$-cube} $GQ(r,n)$, defined on p. 28 of \cite{Kotsis}, is the direct product
$C_r\times H_n$. So, $C_r\times H_n\to \frac{1}{2}H_{r+2n}$ and, for even $r$, $C_r\times H_n\to H_{\frac{r}{2}+n}$.

A {\em Mesh} and a {\em Generalised Hypercube} are direct products of paths and of complete graphs, respectively,
Clearly, it holds
\begin{equation*}
(P_{m_1}\times \dots \times P_{m_k})\to H_{(m_1+ \dots +m_k)-k}\,\, \mbox{~and~} \,\,
(C_{m_1}\times \dots \times C_{m_k})\to \frac{1}{2}H_{(m_1+ \dots +m_k)}.
\end{equation*}

\begin{figure}
\begin{center}
\begin{minipage}{6cm}
\begin{center}
{\small
\begin{tabular}{||c|c||}
\hline
\hline
1  : (0,0,0,0,0,0,0,0) & 2  : (1,1,0,0,0,0,0,0)\\
3  : (1,1,0,0,1,1,0,0) & 4  : (1,1,1,1,1,1,1,1)\\
5  : (1,1,0,0,1,1,1,1) & 6  : (1,1,0,0,0,0,1,1)\\
7  : (1,1,1,1,0,0,0,0) & 8  : 0,0,1,1,0,0,0,0)\\
9  : (0,0,1,1,1,0,1,0) & 10 : (0,0,0,0,1,1,1,1)\\
11 : (0,0,1,1,1,1,1,1) & 12 :  (0,0,1,1,0,1,0,1)\\
\hline
\end{tabular}
}
\end{center}
\end{minipage}
\begin{minipage}{6cm}
\begin{center}
\centering
\resizebox{25mm}{!}{\includegraphics{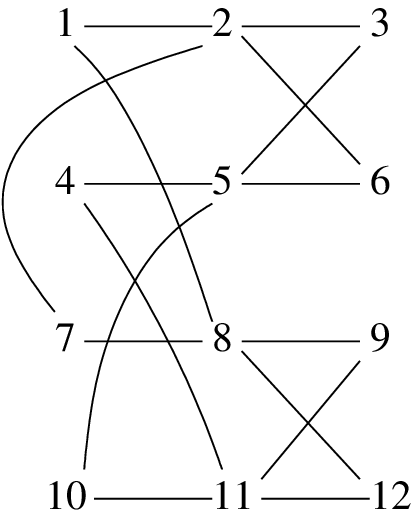}}\par
\end{center}
\end{minipage}
\end{center}
\caption{Butterfly graph $But(2)$ and a $3$-embedding  $But(2)\to \frac{1}{2}H_8$}
\label{tab7}
\end{figure}

The undirected {\em Butterfly Graph} $But(n)$ is (cf., for example, p. 12 in \cite{Wes}) a
 graph on  $2^n(n+1)$ pairs $(x,i)$, where $x$ is a binary sequence of length $n$ and $i\in \{0,1, \dots , n\}$, with
 vertices $(w,i)$ and $(w',i+1)$ being adjacent  if $w'$ is identical to $w$ in all bits with the possible exception of the $(i+1)$-th bit counted from the left. (Note that the definition of Butterfly Graph in \cite{He} is slightly different: it has $2^nn$ vertices there.) The diameter of $But(n)$ is $2n$. It holds
 $But(1)=C_4=H_2$, while $But(2)$ and  $But(3)$ are not $4$-tr. embeddable.
Still $But(2)$ admits nine $3$-tr.embeddings into $\frac{1}{2}H_{8}$;
see one of them on Table \ref{tab7}. Each column of this $12 \times 8$ binary matrix  has exactly $6$ ones.  

The {\em Fibonacci Cube} $Fi(n)$ is the subgraph of $H_n$ induced by the binary {\em Fibonacci sequences}, i.e.,  those containing no two consecutive ones. The {\em Lucas cube} $Lu(n)$ is the subgraph of 
$H_n$ induced by Fibonacci sequences $x_1, \dots , x_n)$ such that not both $x_1$ and $x_n$ are equal to $1$.   Both, $Fi(n)$ and $Lu(n)$, are partial cubes; cf. \cite{Kla}.

\section{Cayley graphs on $Sym(n)$}\label{Cay}
 
 Given  a finite group $G$ and  a generating set $S$ with $S=S^{-1}$ and $id\notin S$,
 the {\em Cayley graph} $CG(G,S)$ having $G$ as the vertex-set and the edge-set consists of pairs of the form $(g, gs)$, with $g\in G, s\in S$. 
 Most of the vertex-transitive  structures - $m$-cubes, Generalised $m$-cubes, Cube-Connected-Cycles (but not Petersen graph) - are Cayley graphs. Below $G$ will always be the symmetric group $Sym(n)$.

The  {\em Star graph} $SG(n)$ is the Cayley graph with $S=\{(1,2),(1,3),\dots, (1,n)\}$;
its diameter is $\left\lfloor \frac{3(n-1)}{2}\right\rfloor$.
It holds $SG(3)=BSG(3)$.
We conjecture that $SG(n)$ with $n\ge 4$ is not $3$-tr.embeddable and checked it for $n=4,5$ and $6$.
 
The {\em Bubble Sort graph} $BSG(n)$ is  the Cayley graph
with $S=\{(1,2),(2,3), \dots ,(n-1,n)\}$; its diameter is  ${n\choose 2}$. 
Its geodesic metric is  called {\em Kendall $\tau$ distance} $I(x,y)$.
 (or {\em inversion metric}, {\em permutation swap metric}).
It is an editing metric on $Sym(n)$:
 the number of adjacent transpositions needed to obtain $x$ from $y$. Also, $I(x,y)$ is the number of {\em relative inversions} of $x$ and $y$, i.e., pairs $(i,j), 1 \le i < j \le n$, with $(x_i-x_j)(y_i -y_j)< 0$. 
\begin{proposition}
For any Bubble Sort graph, it holds  $BSG(n)\to H_{{n\choose 2}}$. 
\end{proposition} 
Proposition 1 in \cite{DDSS} (cf. also Table 3 there) shows that, given a finite Coxeter group $W$ and
its canonical generating set $S$, the Cayley graph $Cay(W, S)$ is isometrically embeddable into $H_{|T|}$,
where $T$ is the set of elements, that are conjugate to an element of $S$.
Above Proposition is just the case $W={\mathsf A}_{n-1}$, since $Sym(n)$ is isomorphic to the the
finite Coxeter group ${\mathsf A}_{n-1}$.
  
\begin{figure}
\begin{center}
\resizebox{50mm}{!}{\includegraphics{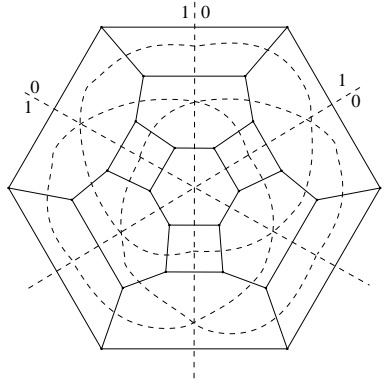}}\par
\end{center}
\caption{Embedding of Bubble Sort graph: $BS(4)\to H_6$} 
\end{figure}

 The {\em Pancake graph} $Pc(n)$ is  the Cayley graph
 with $S$ consisting of $n-1$ permutations of the form $(i,i-1,\dots, 1,i+1,, \dots ,n)$; cf. \cite{He}.
 It holds $$Pc(3)=C_6\to H_3.$$ We conjecture that for $n\ge 4$, $Pc(n)$ is not $3$-tr.embeddable and checked it for the cases $n=4,5,6$
  with diameters $4,5$ and $6$, respectively. To find the diameter of $Pc(n)$ in general, is an open problem, called
  the {\em prefix reversal problem}.  
 
The 
{\em Swap-or-Shift graph} $SOS^n_n$ is the Cayley graph with $S$ consisting the shift $(1, \dots , n)$
and transposition $(1,2)$. The 
graph $SOS^{n-1}_n$ is the Cayley graph with $S$ consisting the shift $(2, \dots , n)$
and transposition $(1,2)$. 
It holds $$SOS^3_3
=Prism_3\to \frac{1}{2}H_5 \mbox{~and~} SOS^4_4\to \frac{1}{2}H_{12},$$ but $SOS^5_5$, having diameter $10$, is not $6$-tr.embeddable.
 It holds $$SOS^2_3=C_6\to H_3,$$ 
but $SOS^3_4$, having diameter $6$, is not $5$-tr.embeddable, but it admits four $4$-tr.embeddings
into $\frac{1}{2}H_{14}$. Also, $SOS^4_5$, having diameter $9$, is not $5$-tr.embeddable. 

So, we expect  that 
$SOS^{n}_n$ is not $6$-tr.embeddable for $n\ge 5$ and 
$SOS^{n-1}_n$ is not $5$-tr.embeddable for $n\ge 4$.

\begin{table}
\begin{center}\small
\begin{center}
\begin{tabular}{||c|c||}
\hline
\hline
ABCD : (0,0,0,0,0,0) & DBCA : (0,0,1,1,1,0)\\
DACB : (1,1,0,1,0,1) & DBAC : (1,0,0,0,0,1)\\
ADCB : (1,1,0,1,1,1) & BDCA : (0,1,1,1,1,0)\\
BACD : (1,0,0,0,0,0) & BDAC : (1,1,0,0,0,1)\\
ADBC : (0,0,1,0,1,0) & CDBA : (1,0,0,0,1,1)\\
CABD : (0,1,1,1,0,1) & CDAB : (0,1,1,0,0,0)\\
ABDC : (0,1,1,1,1,1) & CBDA : (0,1,0,0,0,1)\\
CADB : (1,0,1,0,1,0) & CBAD : (1,0,1,1,1,1)\\
ACDB : (1,0,0,0,1,0) & BCDA : (0,1,0,0,0,0)\\
BADC : (1,1,1,1,1,1) & BCAD : (1,0,1,1,1,0)\\
ACBD : (0,1,0,1,0,1) & DCBA : (1,0,0,1,1,1)\\
DABC : (0,0,1,0,0,0) & DCAB : (0,1,1,1,0,0)\\
\hline
\end{tabular}
\end{center}
\begin{center}
\begin{tabular}{||c|c||}
\hline
\hline
ABCD : (0,0,0,0,0,0,0,0,0,0,0,0,0,0) & DBCA : (1,0,1,0,0,0,0,0,1,0,1,0,1,1)\\
DACB : (0,0,0,0,1,1,0,0,0,1,1,0,1,1) & DBAC : (1,0,1,0,0,0,1,1,1,0,1,0,1,1)\\
ADCB : (1,1,0,0,0,0,1,1,1,0,1,0,0,0) & BDCA : (0,0,0,0,1,1,1,1,0,0,1,1,0,0)\\
BACD : (1,0,0,1,1,1,1,1,1,0,1,0,1,1) & BDAC : (0,0,0,0,1,1,0,0,0,0,1,1,0,0)\\
ADBC : (1,1,0,0,0,0,1,1,1,0,1,0,1,1) & CDBA : (1,0,1,0,0,0,0,0,0,1,1,0,1,1)\\
CABD : (1,1,0,0,0,0,0,0,0,0,0,0,0,0) & CDAB : (1,0,1,0,1,1,0,0,0,1,1,0,1,1)\\
ABDC : (0,0,0,0,1,1,0,0,0,0,0,0,0,0) & CBDA : (1,0,0,1,1,1,1,1,0,0,1,1,0,0)\\
CADB : (1,1,0,0,0,0,1,1,0,0,0,0,0,0) & CBAD : (1,0,0,1,1,1,1,1,0,0,1,1,1,1)\\
ACDB : (1,0,0,1,1,1,0,0,0,1,1,0,1,1) & BCDA : (1,0,1,0,0,0,0,0,0,0,0,0,1,1)\\
BADC : (1,0,0,1,0,0,1,1,1,0,1,0,1,1) & BCAD : (1,0,1,0,0,0,0,0,0,0,0,0,0,0)\\
ACBD : (1,0,0,1,1,1,1,1,0,1,1,0,1,1) & DCBA : (1,1,0,0,1,1,1,1,0,0,1,1,0,0)\\
DABC : (0,0,0,0,1,1,0,0,0,1,1,0,0,0) & DCAB : (1,1,0,0,0,0,1,1,0,0,1,1,0,0)\\
\hline
\end{tabular}
\end{center}
\end{center}
\caption{Embedding  $SOS^4_4\to H_6$ and $4$-tr-Embedding  $SOS^3_4\to \frac{1}{2}H_{14}$}
\label{tab3_SOSab}
\end{table}

\section{Graphs on cycles}\label{Cyc}

The {\em M\"obius ladder} $M_{2m}$ is a cubic circulant graph with $2m$ vertices, formed from an $m$-cycle by adding edges connecting 
opposite pairs of vertices in the cycle. We conjecture that $M_{2m}$ is not $2$-tr.embeddable and checked it  for the {\em Thomsen} (or {\em utility})  {\em graph} $M_6=K_{3,3}$,  Wagner graph $M_8$ (both of diameter $2$) and for $M_{10}$ of diameter $3$.

For  even $n>$ and  increasing sequence $\vec{a}=(a_1, a_2, \dots, a_k)$ of odd numbers from $[3,n-1]$,
we introduce the {\em Generalised Chordal Ring} $GCR(n,\vec{a})$ as the graph obtained by adding to
the cycle $C_{1, \dots, n}$, where  each $i$ is adjacent to $i-1$ and $i+1$ modulo $n$, the following edges:
\begin{enumerate}
\item if $i$ is even, then $i$ is adjacent to $i + a_l \mod n$ for $1\leq l\leq k$;
\item if $i$ is odd, then $i$ is adjacent to $i - a_l \mod n$ for $1\leq l\leq k$.
\end{enumerate}
The cases $k=1$ and $2$ correspond to known topologies: the {\em Chordal Rings} and {\em Double Chordal Rings}, respectively.
The Chordal Ring $GCR(n,a)$ is embeddable for $a=1$ and $3$ (being $C_n$ and $Prism_{\frac{n}{2}}$, respectively), but for $a=5$ and $7$, $GCR(n,a)$ (of diameter $d=3$ and $4$, respectively) is not $(d-1)$-tr.embeddable even for the smallest case $n=2a$.

The results of our computations are summarized in the Conjecture below
and Table \ref{tab2_GCR}, listing known embeddings, which are not covered by this Conjecture (ii).

\begin{figure}
\begin{minipage}{6cm}
\begin{center}
{\small
\begin{tabular}{||c|c||}
\hline
\hline
1 : (0,0,0,0,0) & 2 : (1,0,0,0,0)\\
3 : (1,0,1,0,0) & 4 : (1,0,1,0,1)\\
5 : (1,0,1,1,1) & 6 : (1,1,1,1,1)\\
7 : (1,1,0,1,1) & 8 : (0,1,0,1,1)\\
9 : (0,1,0,0,1) & 10 : (0,1,0,0,0)\\
11 : (0,1,1,0,0) & 12 : (0,0,1,0,0)\\
13 : (0,0,1,1,0) & 14 : (1,0,1,1,0)\\
15 : (1,0,0,1,0) & 16 : (1,0,0,1,1)\\
17 : (1,0,0,0,1) & 18 : (1,1,0,0,1)\\
19 : (1,1,1,0,1) & 20 : (0,1,1,0,1)\\
21 : (0,1,1,1,1) & 22 : (0,1,1,1,0)\\
23 : (0,1,0,1,0) & 24 : (0,0,0,1,0)\\
\hline
\end{tabular}
}
\end{center}
\end{minipage}
\begin{minipage}{6cm}
\resizebox{6cm}{!}{\includegraphics[width=50mm,bb=136 264 489 541,clip]{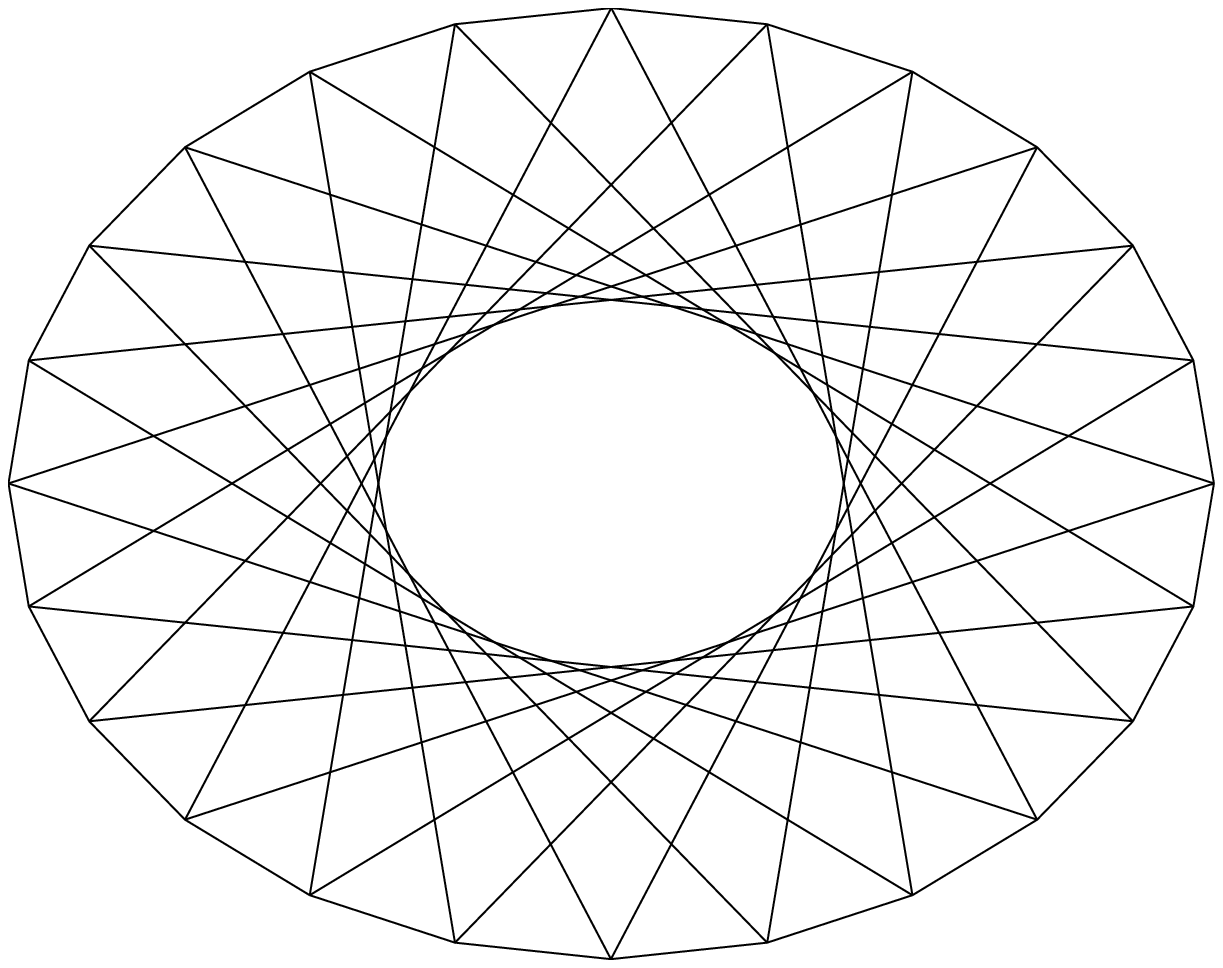}}\par
\end{minipage}
\caption{Embedding of Double Chordal Ring: $GCR(24,(9,11))\to H_5$}
\label{FigChordal}
\end{figure}

See Figure \ref{FigChordal} for $GCR(24, \{9,11\})$, the smallest case in Conjecture (ii).

 \begin{conjecture} (checked for $v\le 70, k\le 5$ and  $v\le 200, k=2, a_2=a_1+2$)
 \label{Conj1}
 
(i) If $GCR(n,\vec{a})$ of diameter $d$ is embeddable, then $n\equiv 0 \mod 4$, $\vec{a}=\{a,a+2\}$
and embedding is into $H_d$.
 
(ii) For each $n\equiv 8 \mod 16, n\ge 24$,
the Double Chordal Rings $GCR(n, (\frac{n}{2}-3, \frac{n}{2}-1))$ and $GCR(n, (\frac{n}{2}+1, \frac{n}{2}+3))$
have $d=\frac{n}{8}+2$ and embed into $H_d$.
 \end{conjecture}

\begin{table}
\begin{center}
{\small
\begin{tabular}{||c|c|c||c|c|c||}
\hline
\hline
n & $\vec{a}$ & emb.  into & n & $\vec{a}$ & emb.  into\\
\hline
48 & (13,15) & $H_{7}$ & 60 & (21,23) & $H_{8}$\\
80 & (17,19) & $H_{9}$ & 84 & (25,27) & $H_{10}$\\
96 & (33,35) & $H_{11}$ & 112 & (29,31) & $H_{11}$\\
120 & (21,23) & $H_{11}$ & 120 & (37,39) & $H_{13}$\\
132 & (45,47) & $H_{14}$ & 140 & (57,59) & $H_{12}$\\
144 & (33,35) & $H_{13}$ & 156 & (49,51) & $H_{16}$\\
160 & (61,63) & $H_{13}$ & 168 & (25,27) & $H_{13}$\\
168 & (57,59) & $H_{17}$ & 176 & (45,47) & $H_{15}$\\
180 & (37,39) & $H_{14}$ & 192 & (61,63) & $H_{19}$\\
\hline
\hline
\end{tabular}
}
\end{center}
\caption{All known embeddings of   $GCR(n,\vec{a})$, not covered by Conjecture \ref{Conj1} (ii). For two $120$- and two $168$-vertex graphs, $\vec{a}$ is 
$(\frac{n}{4}-9, \frac{n}{4}-7), (\frac{n}{4}+7,   \frac{n}{4}+9)$ and $ (\frac{n}{4}-17, \frac{n}{4}-15), (\frac{n}{4}+15,   \frac{n}{4}+17)$}
\label{tab2_GCR}
\end{table}

\section{Regular maps}\label{Map}

A {\em map} is a  {\em $2$-cell decomposition} of a closed compact two-dimensional
 manifold, i.e., a
decomposition of a $2$-manifold  into topological disks. A {\em regular map} is a map
such that every {\em flag} (an incident vertex-edge-face triple) can be transformed into any other flag by a symmetry of the decomposition.
The map {\em of type $\{a,b\}$} is  the regular map with degree $a$ of vertices, having only $b$-gonal faces.

\begin{figure}\begin{center}
\begin{minipage}{6cm}
\centering
\epsfig{file=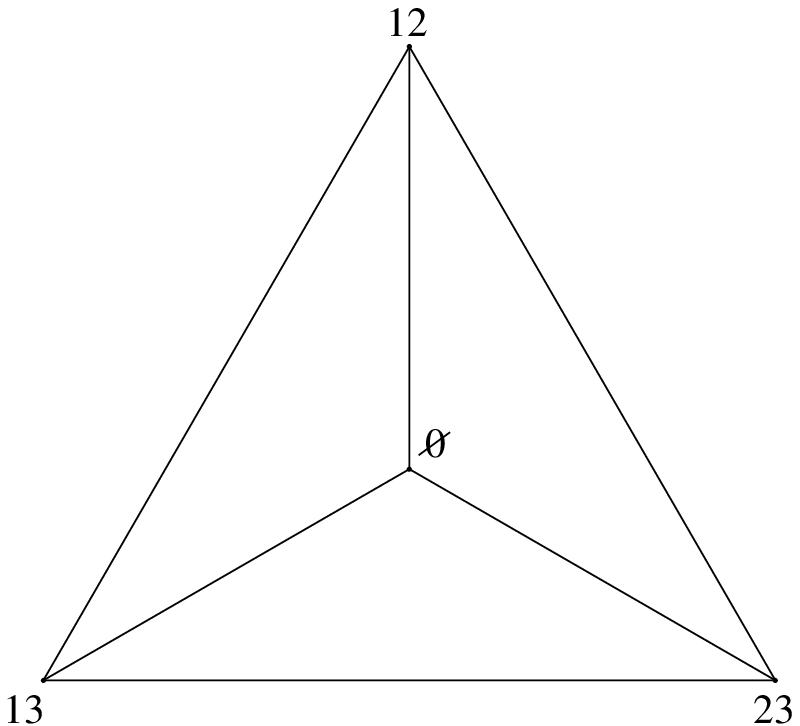, height=3.5cm}\par
\end{minipage}
\begin{minipage}{6cm}
\centering
\epsfig{file=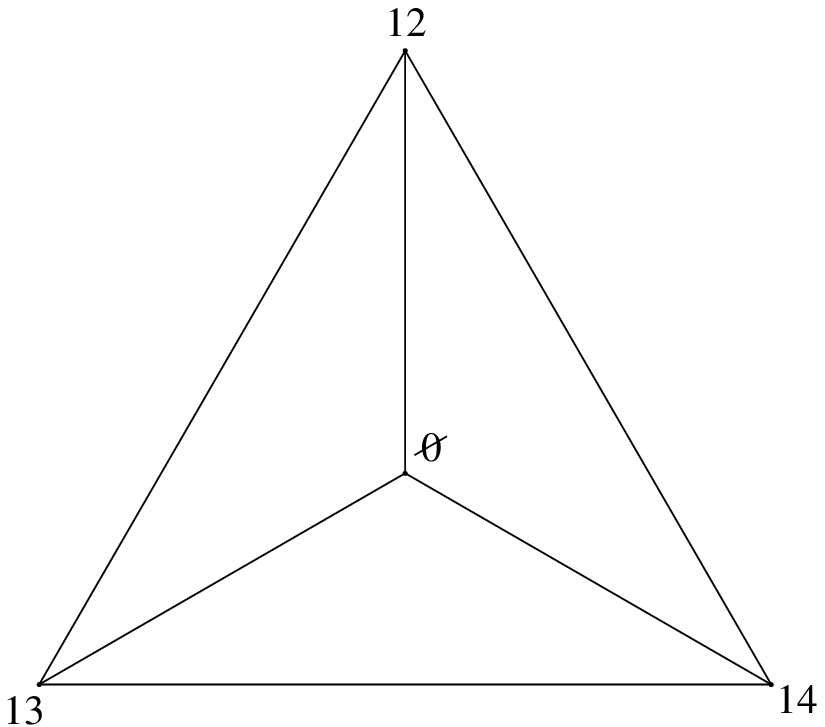, height=3.5cm}\par
\end{minipage}
\begin{minipage}{6cm}
\centering
\epsfig{file=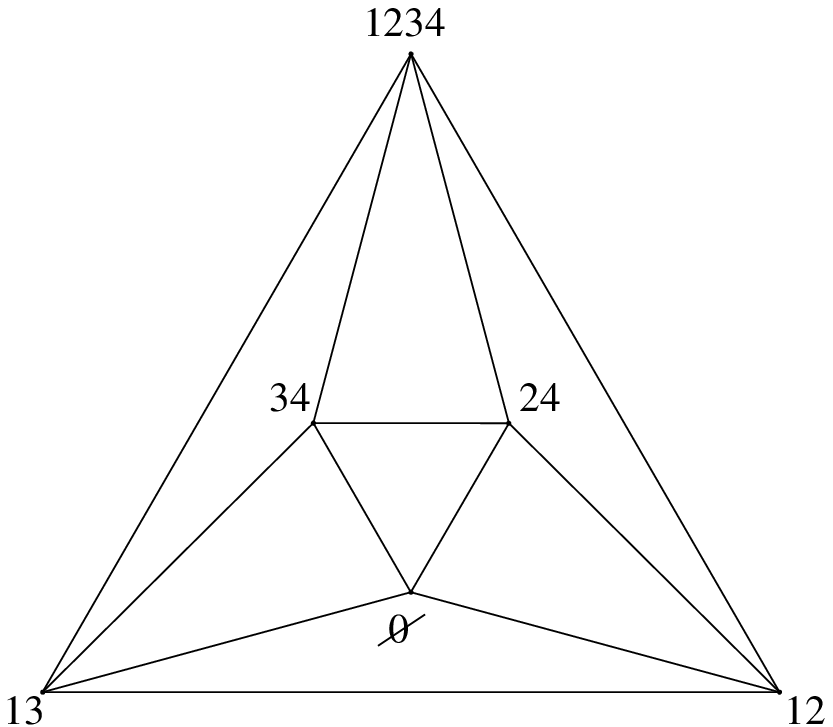, height=3.5cm}\par
\end{minipage}
\begin{minipage}{6cm}
\centering
\epsfig{file=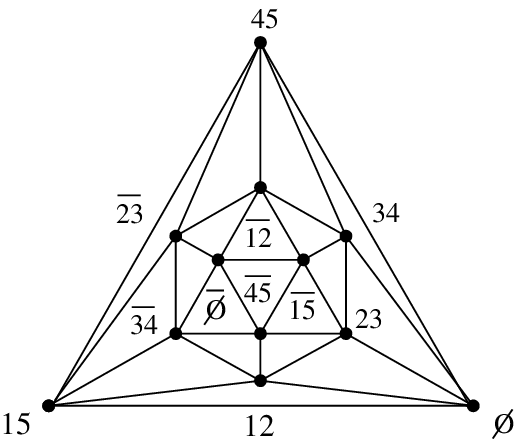, height=3.5cm}\par
\end{minipage}
\end{center}
\caption{Embeddings of Tetrahedron (two), Octahedron and Icosahedron into $3-,4$-, $4$- and $6$-half-cube, respectively}
\label{fig:Embeddings}
\end{figure}

Each of five regular spherical maps, 
 i.e., skeletons of  Platonic polyhedra, are embeddable; cf., say,
  \cite{DGS}.  It holds 
  $$K_4=\frac{1}{2}H_3\simeq J(4,1),\,\,\,\,\,\, K_2^3=H_3,\,\,\,\,\,\,  
    K_{2,2,2}=J(4,2),\,\,\,\,\, \,
K_{2,2,2,2}=\frac{1}{2}H_4$$
for Tetrahedron, Cube, Octahedron, Hyperoctahedron, respectively.
Icosahedron and Dodecahedron
embed into $\frac{1}{2}H_6$, 
 $\frac{1}{2}H_{10}$. See Figs. \ref{fig:Embeddings} 
  and \ref{DurerOctahedronElogatedDodecahedron}.

The {\em cubic Klein graph} is  a $3$-regular graph of diameter $6$ with $56$ vertices, which
 is 
 the skeleton of the {\em Klein map}, 
 a symmetric tessellation of a genus $3$ surface  by $24$ heptagons. 
Neither it, nor its dual are embeddable.
The cubic Klein and Dick graphs are Cayley graphs. See  their dual (on genus $3$ surface) on  Fig. \ref{Dyck}.

 The {\em Dyck graph} is  a $3$-regular graph of diameter $5$ with $32$ vertices, which is 
 the skeleton of the {\em Dyck map}, 
 a symmetric tessellation of a genus-$3$ surface  by $12$ octagons. Neither it, nor
 $K_{4,4,4}$ (its dual for this tiling) are embeddable. 
 The  Dyck graph is toroidal; the skeleton of its dual on the torus $\mathbb{T}^2$ is the {\em Shrikhande graph}, which embeds into $\frac{1}{2}H_6$. The Shrikhande graph can be constructed as a Cayley graph on $\mathbb{Z}_4 \times \mathbb{Z}_4$ with two vertices being adjacent if  the difference is in $\{\pm( 1,0),\pm(0,1),\pm (1,1)\}$.

 The Dyck graph admits a $4$-tr.embedding into $H_6$; see it on Table \ref{tab3}. 
  Each column of this $32\times 6$ binary matrix $X=((x_{ij}))$ has exactly $16$ ones. Let $X'=((1-x_{ij}))$. Clearly, the $32$-sets of rows of $X$ and  $X'$ form together the $64$ vertices of $H_6$. 
 For every vertex  $v$ of the Dick graph, its stabilizer has two orbits, of sizes $3$ and $1$, of antipodal (i.e., at distance $5$) points, say, $\{v'\}$. The distance matrices $D$ of the Dick graph and $D(X)$ (Hamming pairwise distances  of rows of $X$)
  differ only in $16$ entries: $16$ distances of the form $d(v,v')$ are $5$ in $D$, but became $3$
 in $D(X)$.

We analyzed all regular maps from \cite{Con} up to genus $13$ and found embeddings of skeletons for many of them. We do not take just the maps occurring there, but also the maps obtained from  them by the so-called  {\em Wythoff construction} (see \cite{DDSS} for an exposition).
In our context, the Wythoff construction takes a map $M$, a non-trivial subset $S$ of $\{0,1,2\}$ and returns another map $W_{S}(M)$. We embedding

\begin{conjecture}(checked for all maps of genus $g\le 13$)

(i) For any $g\geq 3$, there exist a unique map $M$ of genus $g$ and type $\{4, 4g\}$ such that its skeleton is the cycle $C_{2g}\to H_g$.

(ii) For any $g\geq 2$ there exist a unique map $M$ of genus $g$ and type $\{4g, 4g\}$ such that its skeleton is a cycle $C_{2g}\to H_g$.

(iii) For any $g\geq 2$, there exist a unique map $M$ of genus $g$ and type $\{4, 2g+2\}$ such that its skeleton is a cycle $C_{2g+2}\to H_{g+1}$, the dual skeleton
is $C_4$ and the map $W_{\{0,2\}}(M)$ has  $(8g+8)$-vertex skeleton of diameter $g+3$, that   is embeddable into $H_{g+3}$.

(iv) For any $g\geq 2$, there exist a unique map $M$ of genus $g$ and type $\{2g+1, 4g +2\}$ such that the skeleton of $W_{\{1\}}(M)$ is a cycle  $C_{2g+1}\to \frac{1}{2}H_{2g+1}$
 and the 
 $W_{\{0,1\}}(M)$ has $(4g+2)$-vertex skeleton of  diameter $g+1$, that is embeddable into $\frac{1}{2}H_{2g+3}$.

(v) For any $g \geq 2$, there exist a unique map $M$ of genus $g$ and type $\{2g+2, 2g+2\}$ such that  $W_{\{0,1\}}(M)$ has $(4g+4)$-vertex skeleton 
of diameter $g+2$, that is embeddable into $H_{g+2}$.

\end{conjecture}

\begin{table}
\begin{center}
{\small
\begin{tabular}{||c|c|c|c||}
\hline
\hline
(0,0,0,0,0,0) & (1,0,0,0,0,0) & (1,0,0,1,0,0) & (1,0,0,1,0,1)\\
(1,1,0,1,0,1) & (0,1,0,1,0,1) & (0,1,0,0,0,1) & (0,1,0,0,0,0)\\
(1,0,0,0,1,0) & (1,1,0,0,1,0) & (1,1,1,0,1,0) & (1,1,1,1,1,0)\\
(1,1,1,1,0,0) & (1,0,1,1,0,0) & (0,0,1,0,0,0) & (0,0,1,0,0,1)\\
(0,0,1,0,1,1) & (1,0,1,0,1,1) & (1,0,0,0,1,1) & (0,1,0,0,1,0)\\
(0,1,0,1,1,0) & (0,1,1,1,1,0) & (0,0,1,1,1,0) & (0,0,1,1,0,0)\\
(0,1,1,0,0,1) & (1,1,1,0,0,1) & (1,1,1,0,1,1) & (0,1,0,1,1,1)\\
(0,0,0,1,1,1) & (0,0,1,1,1,1) & (1,1,1,1,0,1) & (1,0,0,1,1,1)\\
\hline
\end{tabular}
}
\end{center}
\caption{$4$-tr.embedding of the Dyck graph into $H_6$}
\label{tab3}
\end{table}

\end{document}